\documentclass{article}

\usepackage{arxiv}
\usepackage{bm}
\usepackage[numbers]{natbib}
\usepackage{algorithm}
\usepackage{url} 
\usepackage{subcaption}
\usepackage{amssymb}
\usepackage{amsmath}
\usepackage{algpseudocode}
\usepackage[utf8]{inputenc}
\usepackage[T1]{fontenc}
\usepackage{hyperref}
\usepackage{url}
\usepackage{booktabs}
\usepackage{amsfonts}
\usepackage{nicefrac}
\usepackage{microtype}
\usepackage{soul}
\usepackage{graphicx}
\usepackage{natbib}
\usepackage{doi}
\usepackage{cleveref}
\usepackage{xcolor}
\usepackage{caption}
\hypersetup{
    colorlinks,
    linkcolor={red!50!black},
    citecolor={blue!50!black},
    urlcolor={blue!80!black}
}
\usepackage{titlesec}
\usepackage{lineno}
\titlespacing\section{0pt}{12pt plus 4pt minus 2pt}{0pt plus 2pt minus 2pt}
\titlespacing\subsection{0pt}{10pt plus 4pt minus 2pt}{0pt plus 2pt minus 2pt}
\titlespacing\subsubsection{0pt}{8pt plus 4pt minus 2pt}{0pt plus 2pt minus 2pt}

\usepackage{enumitem}
\setlist{nosep}
\newcommand*\samethanks[1][\value{footnote}]{\footnotemark[#1]}
\newcommand*{\email}[1]{\texttt{#1}}
\newcommand*{\affaddr}[1]{#1}

\title{ FluTO: Graded Multiscale Fluid Topology Optimization using Neural Networks}

\date{}
\author{
Rahul Kumar Padhy \thanks{Contributed equally}, Aaditya Chandrasekhar \samethanks, Krishnan Suresh \\
\affaddr{Department of Mechanical Engineering, University of Wisconsin-Madison}\\
\email{\{rkpadhy, achandrasek3, ksuresh\}@wisc.edu}
}

\renewcommand{\headeright}{}
\renewcommand{\undertitle}{}
\renewcommand{\shorttitle}{GM-Flow}

\begin{document}
\maketitle

\begin{abstract}

Fluid-flow devices with low dissipation, but high contact area, are of importance in many applications. A well-known strategy to design such devices is multi-scale topology optimization (MTO), where optimal microstructures are designed within each cell of a discretized domain. Unfortunately, MTO is computationally very expensive since one must perform  homogenization of the evolving microstructures, during each step of the homogenization process. As an alternate, we propose here a graded multiscale topology optimization (GMTO) for designing fluid-flow devices. In the proposed method, several pre-selected but size-parameterized and orientable microstructures are used to fill the domain optimally. GMTO significantly reduces the computation while retaining many of the benefits of MTO. 

In particular,  GMTO is implemented here using a neural-network (NN) since: (1) homogenization can be performed off-line, and used by the NN during optimization, (2) it enables continuous switching between microstructures during optimization, (3) the number of design variables and computational effort is independent of number of microstructure used, and, (4) it supports automatic differentiation, thereby eliminating manual sensitivity analysis.  Several numerical results are presented to illustrate the proposed framework.

 \begin{figure}[H]
 	\begin{center}
		\includegraphics[scale=0.9,trim={0 0 0 0},clip]{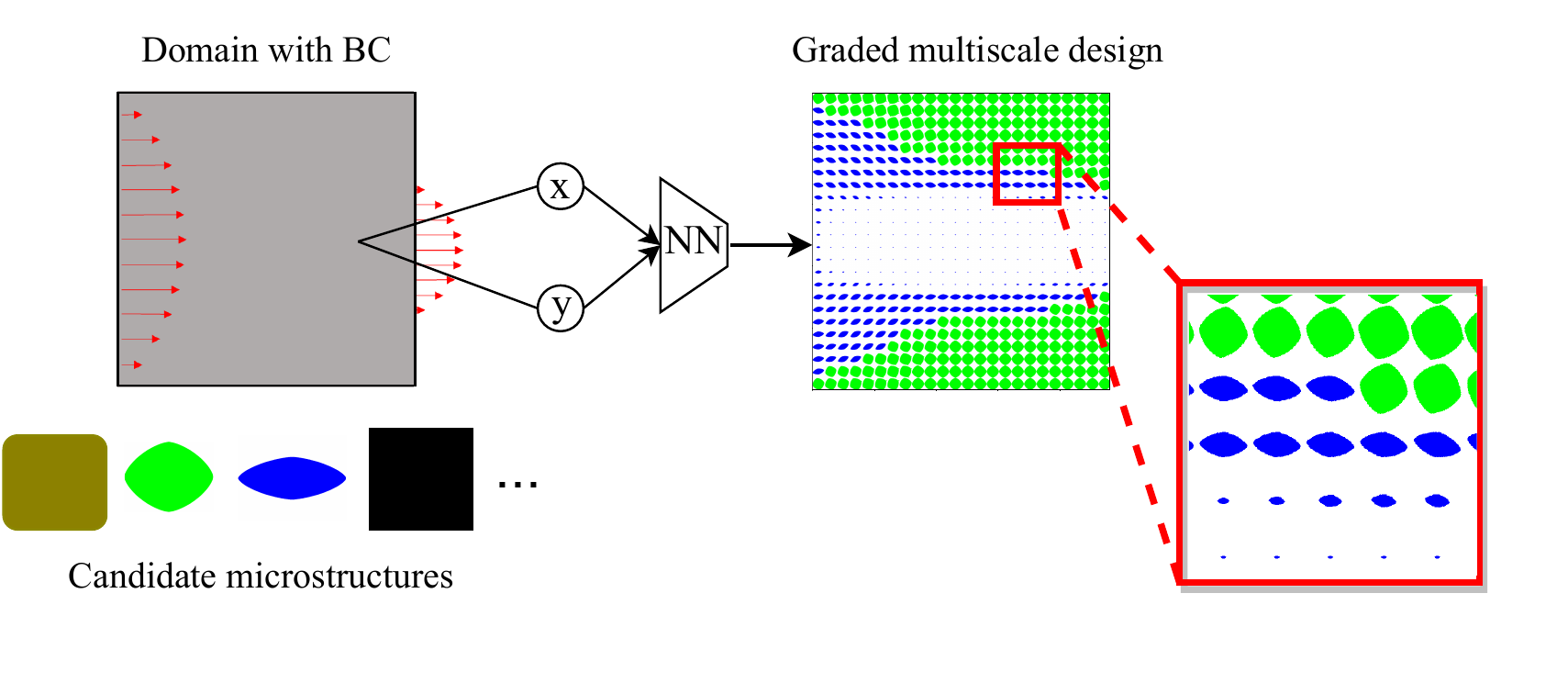}
 		\caption{ Graphical abstract: Given a set of candidate microstructures and a fluid topology optimization problem, a neural network (NN) selects appropriate microstructures, optimizes their size and orientation to produce a graded multiscale design.}
 		\label{fig:abstract}
	\end{center}
 \end{figure}
\end{abstract}


\section{Introduction}
\label{sec:intro}

\begin{figure}[]
 	\begin{center}
		\includegraphics[scale=0.5,trim={0 0 0 0},clip, angle= 270]{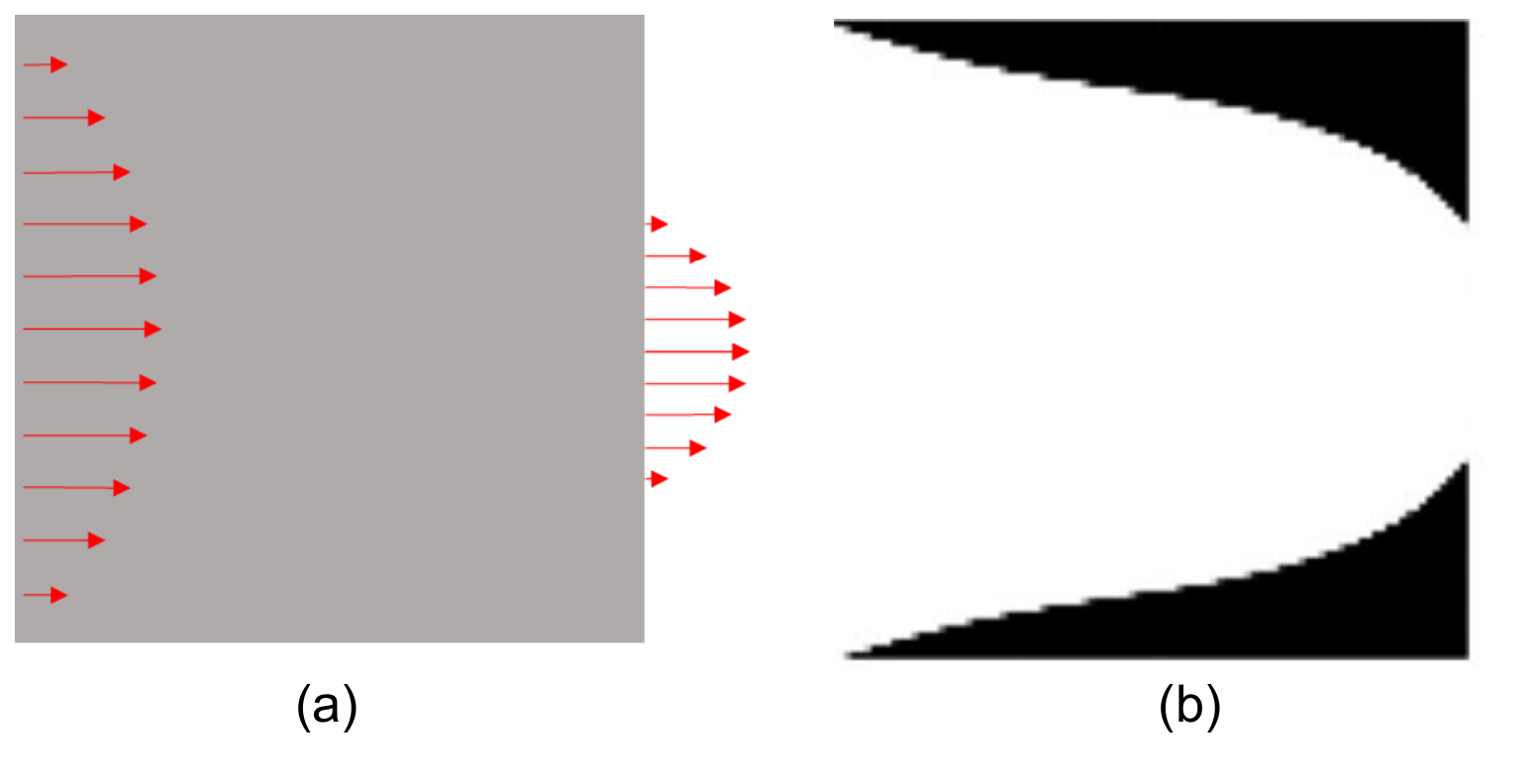}
 		\caption{ Classic fluid-flow topology optimization to minimize dissipated power.}
 		\label{fig:single_scale}
	\end{center}
 \end{figure}

Topology optimization is used extensively today in various disciplines \cite{Sigmund2013}. When applied to  fluid-flow problems, it addresses the following question \cite{alexandersen2020review}: \emph{"Where should  fluid flow within a design domain?"} As an example, consider the simple diffuser problem posed in \cref{fig:single_scale}(a) where fluid enters on the left  and  exits on the right, as illustrated. The objective is to design the  fluid flow path such that the dissipated power is minimized. In addition, a volume constraint is often imposed on the surrounding solid material (or equivalently, on the fluid space). Using topology optimization techniques described, for example in \cite{borrvall2003topology}, one can  solve this problem to arrive at the topology, i.e., fluid flow path, illustrated in \cref{fig:single_scale}(b). For such simple problems, one can argue that topology optimization is not needed since the solution is often a single flow path \cite{alexandersen2020review}.

\begin{figure}[]
 	\begin{center}
		\includegraphics[scale=0.5,trim={0 0 0 0},clip, angle= 270]{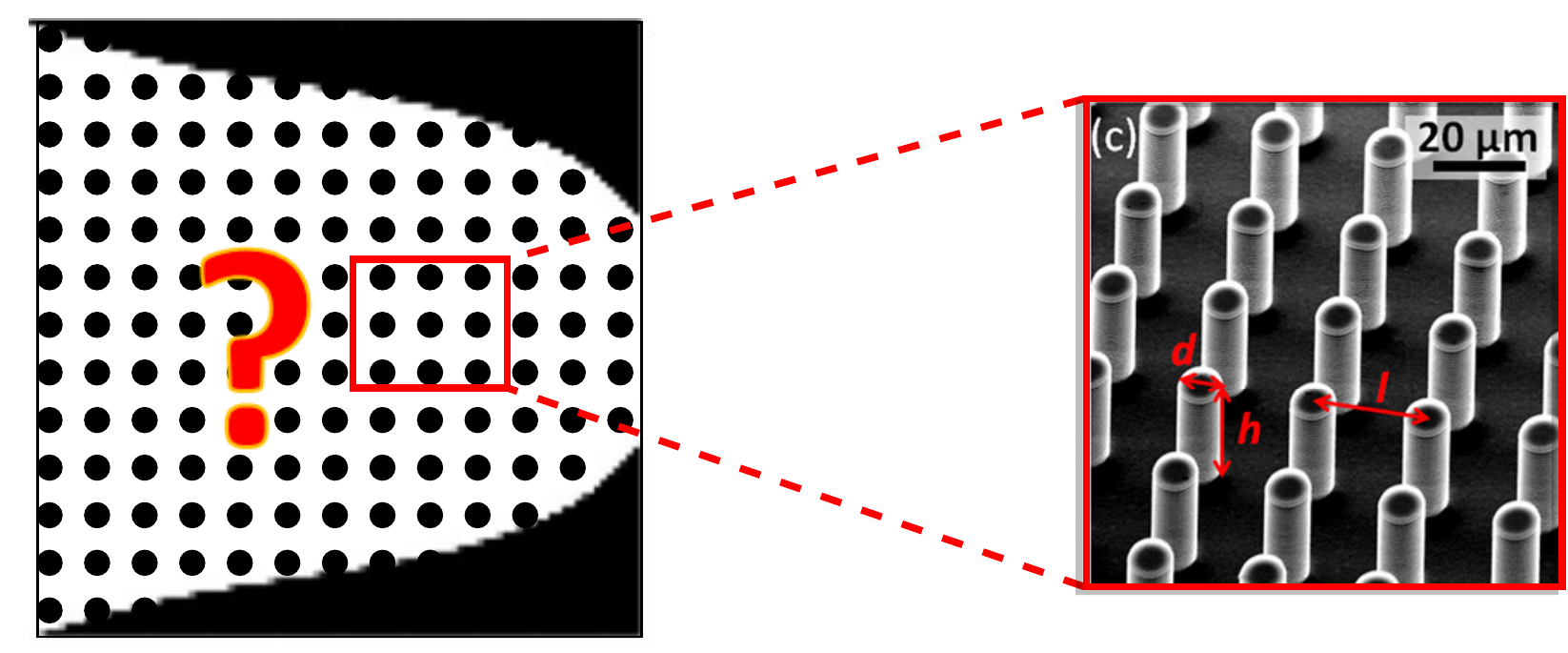}
 		\caption{ Uniform cylindrical micro-pillars \cite{zhu2016prediction}.   }
 		\label{fig:pillars}
	\end{center}
 \end{figure}

Topology optimization becomes more important when additional criteria  come into play. For example, maximizing the fluid-solid  contact area is critical in many  applications including bio-sensors for detecting tumour cells \cite{nagrath2007isolation}, microfluidic devices for cell sorting \cite{fan1999dynamic, hayes2001flow, jiang2000mrna, liu2007micropillar, choi2002integrated}, micro-channel heat sinks \cite{zhu2016prediction, guo2013multiphysics, moran2004microsystem}, and other microfluidic devices involving heat transfer and mass transportation/mixing mechanisms \cite{bixler2012bioinspired,bixler2013fluid}. In such applications, a heuristic strategy to increase the contact surface is to use micro-pillar arrays   \cite{huang2018review, li2014high}  as illustrated in  \cref{fig:pillars}. 

While uniform cylindrical micro-pillars can increase the contact surface, they can also significantly increase the dissipated power loss \cite{lauder2016structure,bocanegra2016holographic}. To balance the two, one can vary the the cylinder radii or, better still, use non-cylindrical micro-pillars, henceforth referred to as \emph{microstructures}, with proper orientation for minimizing power loss and maximizing fluid-solid contact interface \cite{wu2019topology}. Optimizing the cross-section, size and orientation of the microstructures leads to a \emph {multi-scale fluid-flow topology optimization problem} \cite{wu2019topology}, the main focus of this paper.
In multi-scale topology optimization (MTO), optimal microstructures are designed in each cell of a discretized domain. Unfortunately, this can be very expensive \cite{wu2021topology} since one must carry out  homogenization of the evolving microstructures \cite{zhou2008design} during each step of the global optimization process. Specifically, if $N_e$ is the number of elements (typically in the order of 1000s), $N$ is the number of global optimization steps (typically in the order of 100s), and $c$ is the cost of homogenization, then the MTO cost is at least $N * N_e * c$. In this paper, we propose a graded MTO (GMTO) strategy that retains much of the advantages of MTO, but exhibits a significantly lower computational cost.

The remainder of this paper is organized as follows. In \cref{sec:RelWork} we review the literature on topology optimization for fluid flow problems, and summarize the main contributions of this paper. In \cref{sec:tecBack}, we briefly review relevant technical background, followed by  the proposed method in  \cref{sec:method}. Numerical examples are illustrated in \cref{sec:experiments}, and conclusions are provided in \cref{sec:conclusion}.
 
\section{Related Work}
\label{sec:RelWork}

In this section, we briefly review prior work on fluid-flow topology optimization methods. In particular, while many different strategies have been proposed \cite{alexandersen2020review, challis2009level}, we focus here on density-based methods.

\subsection{Single-scale TO}
\label{sec: Min_drag}
Topology optimization for fluid-flow problems first appeared in 2003 in the seminal work of Borrvall and Petersson \cite{borrvall2003topology}. They presented the optimal layout of channel flows for minimal drag (or pressure drop), for Stokes equation with Brinkman–Darcy law equations under low Reynolds number (laminar flow conditions). Gersborg-Hansen et al. \cite{gersborg2005topology} continued this study and presented applications with low Reynolds numbers for microfluidic  problems and micro-electro-mechanical devices. Guest and Prévost \cite{guest2006topology} solved the formulation of the Stokes–Darcy problem numerically using stabilized finite element methods. Wiker et al. \cite{wiker2007topology} used the viscosity as a dependent parameter and presented examples of channels in a tree-shaped structure for a pure Darcy problem and mixed Stokes–Darcy flow. Pereira et al. \cite{pereira2016fluid} reproduced the classical examples presented by \cite{borrvall2003topology} for optimal channels designs considering Stokes–Darcy flow using polygonal meshes and provided an educational software written in MATLAB. Suárez et al.\cite{suarez2022virtual} applied topology optimization to non-Newtonian flows in arbitrary domains, using a virtual element method.

\subsection{Multi-scale TO}
\label{sec: Mstr_db}
Multiscale topology optimization (MTO) for fluid flow, as mentioned earlier, involves generating appropriate microstructures, i.e., generalization of micro-pillars, typically, in each finite element cell, for maximizing fluid permeability. Guest and Prévost \cite{guest2007design} maximised the permeability of porous microstructures using a Darcy–Stokes interpolation \cite{guest2006topology} subject to isotropic symmetry constraints. This was extended in \cite{guest2006optimizing} to optimize microstructures for maximal stiffness and permeability. Bio-mimicking techniques have been demonstrated for achieving this goal \cite{jakvsic2020biomimetic} but are not proven to be optimal \cite{wu2019topology}. MTO that minimizes energy loss and offers high contact area was demonstrated in \cite{wu2019topology}. However, the volume fraction of the microstructures was pre-determined, and the surface contact area was not explicitly controlled. 

To address the high computational cost  of MTO, graded MTO (GMTO) techniques have been proposed, in  structural mechanics \cite{chandrasekhar2022gm, nguyen2021multiscale, zhao2022stress, zheng2021data, wang2021data, wang2022data,watts2019simple, white2019multiscale,wang2017multiscale}, but have not been extended to fluid flow problems.

\subsection{Paper Contributions}
\label{sec: PaperContributions}

\textit{In this work we extend a particular GMTO framework proposed in \cite{chandrasekhar2022gm} for structural problems, to fluid flow problems}. In structural GMTO, a total volume constraint is imposed, whereas here, a total contact area constraint is more critical, changing the type of microstructures one must choose. Further, microstructure orientation is not typically considered in structural problems (due to potential loss in connectivity), but it plays an important role in fluid problems. 

The proposed GMTO strategy is based on a neural-network based optimization framework for the the following reasons: (1) it implicitly guarantees the partition of unity, i.e., ensures that the net volume fraction of microstructures in each cell is unity, as described later on, (2) it supports automatic differentiation, and (3)  the number of design variables is only weakly dependent on the number of pre-selected microstructures.

\section{Technical Background}
\label{sec:tecBack}

\subsection{Fluid Flow Governing Equations}
\label{sec:Fluid_eqs}
We  assume here a low-Reynolds in-compressible Stokes flow, i.e.,
\begin{subequations}
\begin{align}
    -2\nabla.[\mu\bm{\epsilon(u)}]+\bm{C^{-1}}\bm{u}+\nabla p= 0 \text{ in  $\Omega$} \\
    \nabla.\bm{u}=0  \text{ in  $\Omega$}\\
    \bm{u}=\bm{g} \text{ over }  \partial\Omega
    \end{align}
    \label{eq:fluidGovnEq}
\end{subequations}
where $\bm{u}$ and p are the velocity and pressure of the fluid, $\mu$ is its viscosity,  $\bm{\epsilon(u)} = (\bm{\nabla} \bm{u} + \bm{\nabla^T}\bm{u})/2$ denotes the rate-of-strain tensor, $\bm{C^{-1}}$ is the inverse permeability tensor; the fluid is assumed to be of unit density.

\subsection{Problem Formulation}
\label{sec:Problem_specification}

\begin{figure}[H]
 	\begin{center}
	\includegraphics[width=120mm,scale=0.6,trim={0 0 0 0} ]{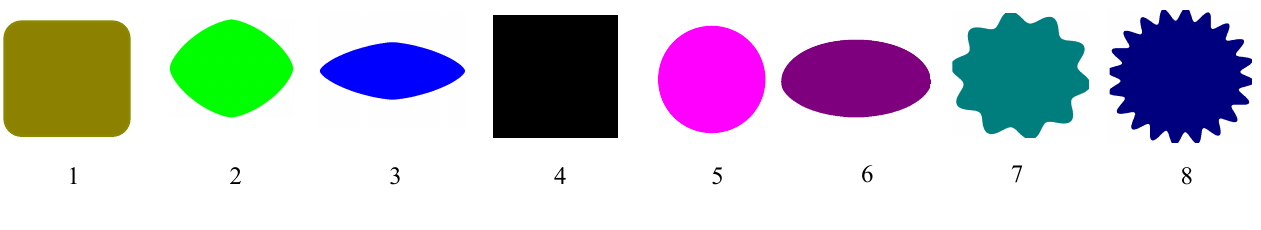}	\caption{Pre-selected microstructures.}
	\label{fig:microstructures}
	\end{center}
 \end{figure}

We also assume that a set of microstuctures have been pre-selected. Specifically, we consider the microstructures illustrated in  \cref{fig:microstructures} (selected from various sources \cite{guest2007design, guest2006topology,aziz2008online,mohammed2016theory, wu2019topology, balbi2013morpho}) that exhibit a wide range of permeability and contact area. Each of these microstructures can be scaled and oriented for optimal flow. For future reference, they are named as follows: 1: Squircle, 2: Fish-body-1, 3: Fish-body-2, 4: Square, 5: Circle, 6: Ellipse, 7: Mucosa-10, and 8: Mucosa-20.  Such microstructures are often seen in nature; for example, mucosa-like structures are known to be present in the human intestine \cite{balbi2013morpho}.

\begin{figure}[]
 	\begin{center}
		\includegraphics[scale=0.5,trim={0 0 0 0},clip]{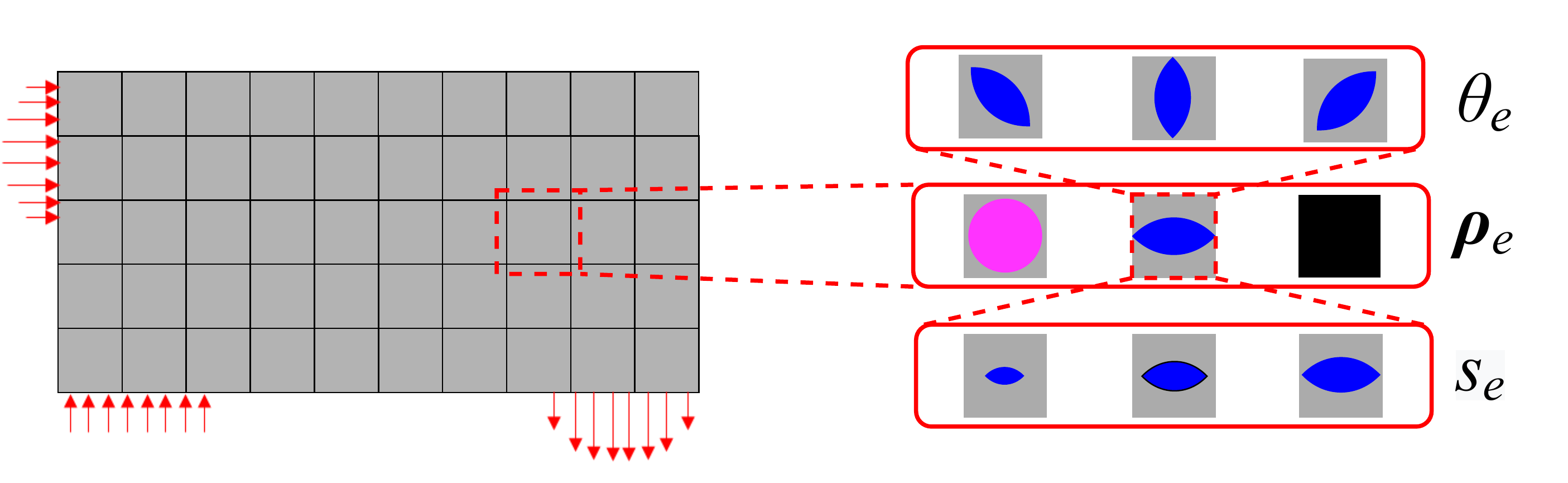}
 		\caption{Graded multiscale TO.}
 		\label{fig:fluidGMTO}
	\end{center}
 \end{figure}

A typical design domain with prescribed flow boundary conditions is illustrated in \cref{fig:fluidGMTO}. We seek to compute an optimal multiscale design where we determine, in each finite element cell, the appropriate microstructure, its size (gradation) and orientation. The objective is to minimize the dissipated power subject to a total contact area (i.e., perimeter in 2D) constraint. The  strategy we adopt is to pre-compute the permeability and contact area of each of the microstructures as a function of its size and orientation, and exploit these for a global multi-scale fluid-flow optimization. 

Towards this end, we introduce the following design variables. The presence or absence of a microstructure at any element $e$ will be captured by the set of density variables  $\bm{\rho}_e = \{\rho_{e,1}, \rho_{e,2}, ..., \rho_{e,M} \}$; see \cref{fig:fluidGMTO}. Ideally $ \rho_{e,m}$ should take a binary value $0/1$ subject to \emph{ partition of unity} constraint  $\sum\rho_{e,m} = 1, \forall e$. However, for gradient-based optimization, we will let $0 \leq \rho_{e,m} \leq 1$ and drive it towards $0/1$ through penalization (once again subject to the partition of unity constraint).

We control the size of all the microstructures at any element $e$  by the scalar design variable $0 \leq s_e \leq 1$ (see \cref{fig:fluidGMTO}), where 0 denotes completely filled with fluid, and 1 denotes maximum size of the microstructure(s). Note that: (1) a single variable $s_e$ controls the size of all microstructures at that point,  (2) $s_e=1$ means that the microstructures at that element are at their maximum size (that can fit into that cell); however, it does not imply that the cell is fully solid, and (3) as we drive $\rho_{m,e}$ towards $0/1$ only one microstructure of that size will prevail at that element.

Finally,  the orientation of all the microstructures at $e$, with respect to the $x$ axis,  will be denoted by $0\leq\theta_e\leq2\pi$ (see \cref{fig:fluidGMTO}). Thus the design variables associated with each element will be denoted by
$\bm{\zeta}_e = \{\rho_{e,1}, \rho_{e,2}, \ldots, \rho_{e,M}, s_e, \theta_e \}$.

\subsection{ Effective Permeability through Numerical Homogenization}
\label{sec:Homogenization}
\begin{figure}[H]
 	\begin{center}		\includegraphics[scale=0.6,trim={0 0 0 0},clip]{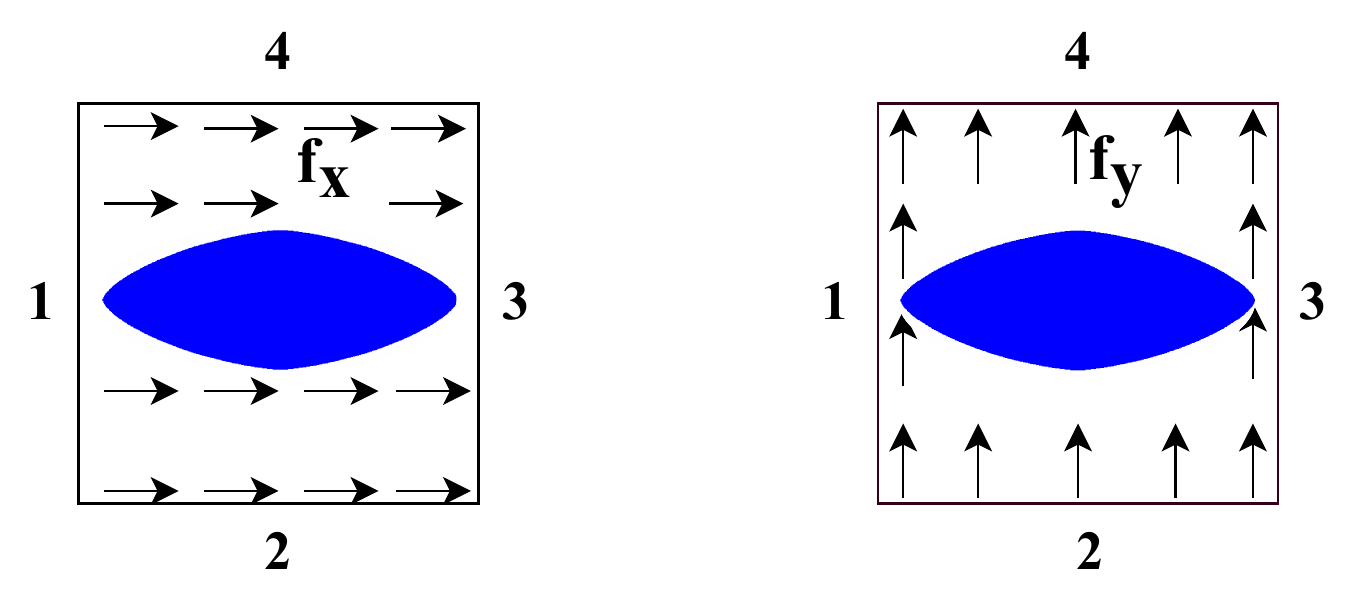}
 		\caption{Fluid flow problems to determine permeability tensor (a) a body force $f_{x} = 1$ to determine velocities {$\{u_{0}(x,y), v_{0}(x,y)\}$} (b) a body force $f_{y} = 1$ to  determine velocities {$\{u_{1}(x,y), v_{1}(x,y)\}$}.}
 		\label{fig:homogenization}
	\end{center}
 \end{figure}

Prior to carrying out global optimization, we pre-compute the  $2 \times 2 $ permeability tensor $\bm{C}_m$ of each microstructure at discrete sizes. The components of $\bm{C}_m$ can be computed by posing two low-Reynolds in-compressible Stokes flow problems over a unit cell, with unit body forces $f_x = 1$ and $f_y = 1$, respectively, as  illustrated in \cref{fig:homogenization}. The boundary conditions for both problems are as follows: (a) boundaries 1 and 3  are coupled through periodic boundary conditions for velocity and pressure, and (b)  boundaries 2 and 4 are similarly coupled. The velocities obtained by solving the problem in \cref{fig:homogenization}a are denoted by $u_0(x,y)$ and $v_0(x,y)$, while those obtained from \cref{fig:homogenization}b are denoted by $u_1(x,y)$ and $v_1(x,y)$.  The components of the permeability tensor $\bm{C}_m$ are then defined as \cite{andreassen2014determine, lang2014permeability, vianna2020computing}: 

\begin{equation}
    \bm{C}_m =  
    \begin{bmatrix} C_m^{00} & C_m^{01} \\
    C_m^{10}  & C_m^{11}
    \end{bmatrix} = 
    \frac{1}{|V|} \begin{bmatrix} \int\limits_{V}u_{0}dV &   \int\limits_{V}v_{0}dV \\
     \int\limits_{V}u_{1}dV &  \int\limits_{V}v_{1}dV \end{bmatrix}
    \label{eq:CMatrix}
\end{equation}
where $V$ is the volume of the unit cell. Observe that, at a default orientation, all the microstructures in \cref{fig:microstructures} are symmetric about the two axis. Therefore, the off-diagonals are zero, i.e., $C_m^{01} =  C_m^{10} = 0$. 

\begin{figure}[H]
 	\begin{center}
		\includegraphics[scale=1,trim={0 0 0 0},clip]{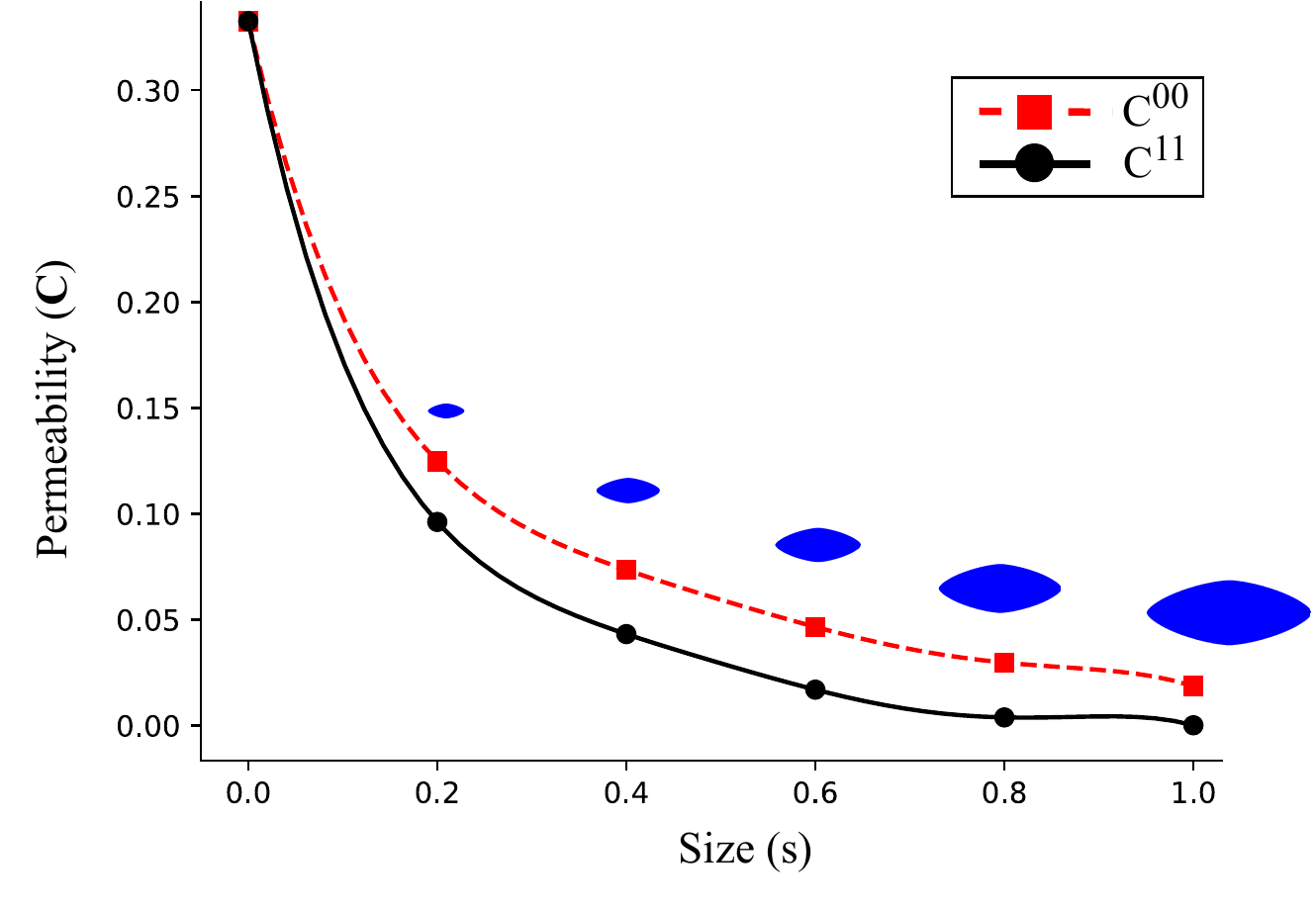}
 		\caption{Permeability components versus size of the fish-body-2 microstructure.}
 		\label{fig:permVsSizeHalfFishscale}
	\end{center}
 \end{figure}

We compute the two remaining components for each microstructure at a finite number of sizes, at default orientation, using the implementation provided in \cite{andreassen2014determine}. Then, two interpolating and positive polynomials \cite{kumar2021spectral} $C_m^{00}(s)$  and $C_m^{11}(s)$  are  constructed from these samples. As an example, the polynomials for a fish-body-type microstructure are illustrated in \cref{fig:permVsSizeHalfFishscale}. In this paper, we use a $5^{th}$ degree polynomial. 

The orientation is then accounted for via the following tensor operation  \cite{lang2014permeability}:
\begin{equation}
    \bm{C}_m(s,\theta) =  \begin{bmatrix} \cos(\theta) & -\sin(\theta) \\ \sin(\theta) & \cos(\theta) \end{bmatrix}
    \begin{bmatrix} C_m^{00}(s) & 0 \\
    0 & C_m^{11}(s) \end{bmatrix}
    \begin{bmatrix} \cos(\theta) & -\sin(\theta) \\ \sin(\theta) & \cos(\theta) \end{bmatrix}^T
    \label{eq:effectiveCMatrix}
\end{equation}
Finally, since multiple microstructure can co-exist at any element, the following penalization scheme is proposed to ensure that we drive $\rho_m$ towards $0/1$ (subscript $e$ has been suppressed for simplicity):
\begin{equation}
    \bm{C}(\bm{\rho}, s, \theta) = \sum\limits_{m=1}^M \rho_m^p \bm{C}_m(s,\theta)
    \label{eq:effectiveCMatrix}
\end{equation}
where the penalization $p > 1$ discourages microstructure mixing. 
 
Observe that the pre-computation cost is approximately $ S * M *c$, where $S$ is the number of size samples (6 in this paper)  , $M$ is the number of microstructures  (8 in this paper) and $c$ is the cost of homogenization. This is significantly lower than the cost for MTO.

\subsection{Contact area }
\label{sec:contact area}
Observe that the contact area $\Gamma_m$ of a microstructure is proportional to its size $s$, and independent of orientation $\theta$. Thus, it is sufficient to compute the maximum contact area $\Gamma_m^{max}$ at $s =1$ for each microstructure. The contact areas are then combined as follows (subscript $e$ has been suppressed for simplicity): 
\begin{equation}
    \Gamma(\bm{\rho}, s) = \sum\limits_{m=1}^M \rho_m s \Gamma_m^{max}
    \label{eq:effectivecontact area}
\end{equation}
\subsection{Volume Constraint}
\label{sec:area_cons}
One can also impose a constraint on the total fluid volume  allowed; this constraint is mainly used in validation experiments (please see \cref{sec:double_pipe_validation}). Similar to the contact area computation, the total volume occupied by the fluid is given by:
\begin{equation}
    V(\bm{\rho}, s) = V_e(1 - \sum\limits_{m=1}^M \rho_m s^{2} v_m^{max})  
    \label{eq:effectivecontact area}
\end{equation}
where $v_m^{max}$ is the fraction of the volume occupied by  a microstructure $m$ at maximum size. 

\subsection{Fluid Flow Finite Element Analysis}
\label{sec:Fluid_FEA}
 
For the global fluid flow analysis, we use  quadrilateral Q2-Q1 (Taylor-Hood) elements. The elemental stiffness matrix $\bm{K_e}$ and degrees of freedom vector $\bm{S_e}$ for the governing equation (see \cref{sec:Fluid_eqs}) are given by (see \cite{pereira2016fluid} for details): 
\begin{equation}
\bm{K_e} = 
\begin{bmatrix}
    \bm{A_e} && \bm{B_e} && \bm{0} \\
    \bm{B_e^T}  && \bm{0} && \bm{a_e}\\
    \bm{0}  && \bm{a_e^T} && \bm{0}
    \end{bmatrix} \;, \; \bm{S_e}=
\begin{bmatrix}
    \bm{U_e}  \\
    \bm{P_e}  \\
    \bm{\lambda}
    \end{bmatrix} 
    \label{eq:StiffnessMatrix}
\end{equation}
where
\begin{subequations}
	\label{eq:stiffness_term_defn}
	   \begin{align}
	& \bm{A_e} = \bm{A_e^\mu}+\bm{C_e}^{-1}\bm{A_e^\alpha}\\
	&\bm{[A_e^\mu]}_{ij} = \int_{\Omega_e}2\mu\bm{\epsilon(N_i)}:\bm{\epsilon(N_j)}d\Omega\\
     &\bm{[A_e^\alpha]}_{ij} = \int_{\Omega_e}\bm{N_i}\bm{N_j}d\Omega\\
     &\bm{[B_e]}_{ij} = \int_{\Omega_e}\bm{L_j}\nabla.\bm{ N_i}d\Omega\\
     &\bm{[a_e]}_{i} = \int_{\Omega_e}\bm{L_i}d\Omega
	\end{align}
\end{subequations}
where $\bm{N_i}$ and $\bm{L_i}$ are the velocity and pressure basis functions, \bm{$U_e$} and \bm{$P_e$} represent elemental velocity and pressure degrees of freedom respectively and $\bm{C_e}  $ is the element permeability matrix (described previously).  In order to uniquely define the pressure field, a zero mean condition is enforced. 

\subsection{Optimization Problem}
\label{sec:Opt}
Consequently, one can pose the GMTO problem in a finite-element setting as:
\begin{subequations}
	\label{eq:optimization_base_Eqn}
	\begin{align}
		& \underset{\overline {\bm{\zeta}} = \{ \bm{\zeta}_1, \bm{\zeta}_2, \ldots \bm{\zeta}_{N_e} \}} {\text{minimize}}
		& J(\overline {\bm{\zeta}}) &= \sum\limits_{e=1}^{N_e}\frac{1}{2}\bm{U_e}^T\bm{[A_e^\mu+A_e^\alpha C_e^{-1}]U_e} \label{eq:optimization_base_objective}\\
		& \text{subject to}
		&  \bm{K}(\overline {\bm{\zeta}})\bm{S} & = \bm{f}\label{eq:optimization_base_govnEq}\\
		& &  g_{\Gamma} (\overline {\bm{\zeta}}) & \equiv 1 - \frac{\sum\limits_{e=1}^{N_e} \sum\limits_{m=1}^M \rho_{e,m} s \Gamma_m^{max}}{\Gamma^*}  \leq 0  \label{eq:optimization_base_perimCons} \\
      & & \text{(or)} \quad  g_{V} (\overline {\bm{\zeta}}) & \equiv   \frac{\sum\limits_{e=1}^{N_e}V_e( 1 -\sum\limits_{m=1}^M \rho_{e,m} v_m^{max}s^{2})}{(\sum\limits_{e=1}^{N_e} V_e) v^*}-1  \leq 0  \label{eq:optimization_base_volCons} \\
		& & \sum\limits_{m=1}^M \rho_{e,m} &= 1 \; , \; \forall e \label{eq:optimization_base_partitionUnity}\\
		& &  0 \leq \rho_{e,m} &\leq 1 \; , \; \forall e \; , \; \forall m \label{optimization_base_boundConsRho} \\
		& &  0 \leq  s_e &\leq 1 \; , \; \forall e \label{optimization_base_boundConsSize} \\
		& & 0 \leq  \theta_e & \leq 2 \pi \; , \; \forall e \label{optimization_base_boundConsTheta}
	\end{align}
\end{subequations}
where $J$ is the dissipated power, ${\Gamma^*}$ is the lower bound on the total contact area, $v^{*}$ is the upper bound on the total volume fraction, $V_e$ is the total unit cell volume, and $\bm{K}(\bm{\zeta})$, $\bm{S}$ and  $\bm{f}$ are the global stiffness, degrees of freedom and the boundary conditions respectively. In the optimization problem, either we impose volume constraint, \cref{eq:optimization_base_volCons}, for validation (please see \cref{sec:double_pipe_validation}) or contact area constraint, \cref{eq:optimization_base_perimCons}, for other numerical experiments.

In the above direct formulation,  the number of design variables is proportional to the mesh size. Further, the partition of unity, and bound constraints must be strictly enforced over each element.  In this paper, we avoid these  issues by indirectly controlling the design variables via a coordinate-based neural-network \cite{chandrasekhar2021tounn,chandrasekhar2021multi}  described next. 


\section{Proposed Method}
\label{sec:method}

\subsection{Neural Network}
\label{sec:method_NN}

\begin{figure}[H]
 	\begin{center}
		\includegraphics[scale=1,trim={0 0 0 0},clip]{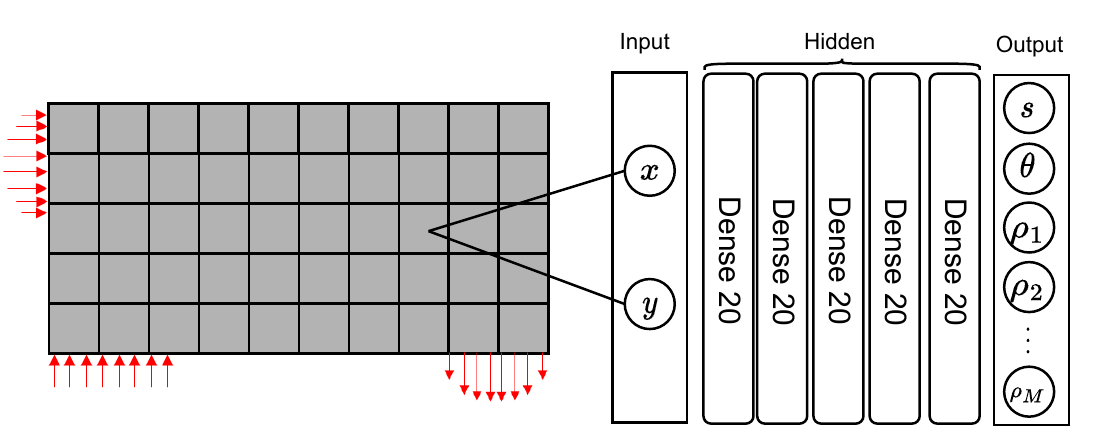}
 		\caption{Neural network architecture.}
 		\label{fig:neuralNetwork}
	\end{center}
 \end{figure}

The proposed neural-network (NN) architecture is illustrated in \cref{fig:neuralNetwork}, and it consists of the following entities:

\begin{enumerate}
	\item \textbf{Input Layer}: The input to the NN are points $\bm{x} \in \mathbf{R}^2$ within the domain. Although these points can be arbitrary,  they correspond here to the center of the elements.
	
	\item \textbf{Hidden Layers}: The hidden layers consist of a series of dense fully connected LeakyReLU activated neurons.
	
	\item \textbf{Output Layer}: The output layer consists of $M + 2$ neurons correspond to the design variables for each element $\zeta = \{\rho_1(\bm{x}), \rho_2(\bm{x}), \ldots, \rho_M(\bm{x}), s(\bm{x}), \theta(\bm{x})\}$. Further, the neurons associated with the density variables are controlled by a softmax function such that the partition of unity $(\sum \rho_m = 1)$ and physical validity $0 \leq \rho_m \leq 1$ are automatically satisfied. The output neuron associated with  the size parameter is controlled via a Sigmoid function $\sigma(\cdot)$, ensuring that $0 \leq s \leq 1$. Finally, the output neuron associated with the orientation parameter is also controlled via a Sigmoid function and scaled as $\theta \leftarrow 2\pi \sigma (\theta)$. Thus, no additional box constraints are needed.
	
	\item \textbf{NN Design Variables}: The weights and bias associated with the NN, denoted by the $\bm{w}$, now become the primary design variables, i.e., we have  $\bm{\rho}(\bm{x}; \bm{w}),$ $s(\bm{x}; \bm{w})$ and $\theta(\bm{x}; \bm{w})$. 
\end{enumerate}

Thus the strategy is to perform GMTO via the NN weights $\bm{w}$, i.e.,  \cref{eq:optimization_base_Eqn} reduces to:
\begin{subequations}
	\label{eq:optimization_nn_Eqn}
	\begin{align}
		& \underset{\bm{w}}{\text{minimize}}
		& &J(\bm{w}) \label{eq:optimization_nn_objective}\\
		& \text{subject to}
		& & \bm{K}(\bm{w})\bm{S} = \bm{f}\label{eq:optimization_nn_govnEq}\\
		& & & g_{\Gamma} (\bm{w}) \equiv  1 - \frac{\sum\limits_{e=1}^{N_e} \sum\limits_{m=1}^M \rho_{e,m} (\bm{w}) s (\bm{w}) \Gamma_m^{max}}{\Gamma^*}  \leq 0  \label{eq:optimization_nn_perimCons}\\
            & &  \text{(or)} \quad  g_{V} (\bm{w}) & \equiv   \frac{\sum\limits_{e=1}^{N_e} V_e(1 -\sum\limits_{m=1}^M \rho_{e,m} (\bm{w}) v_m^{max} s^{2} (\bm{w}))}{(\sum\limits_{e=1}^{N_e} V_e) v^*}-1  \leq 0 \label{eq:optimization_nn_volCons}
	\end{align}
\end{subequations}

Observe that: (1) no additional constraint is needed since they are automatically satisfied, and (2) the number of design variables ($\bm{w}$) is independent of the mesh size and the number of microstructures.

\subsection{Loss Function}
\label{sec:method_loss}
\begin{figure}[H]
 	\begin{center}
		\includegraphics[scale=1,trim={0 0 0 0},clip]{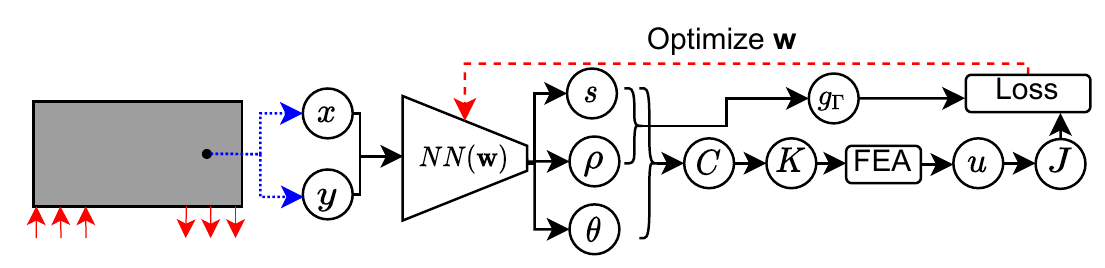}
 		\caption{Optimization loop of the proposed framework.}
 		\label{fig:algoFlowChart}
	\end{center}
 \end{figure}

We now consider solving the NN-based optimization problem in \cref{eq:optimization_nn_Eqn}. Since neural networks are designed to minimize an unconstrained loss function, we  convert the constrained minimization problem into a loss function minimization by employing the augmented Lagrangian scheme  \cite{bertsekas2014constrained}. Specifically, the loss function is defined as
\begin{equation}
    L(\bm{w}) = \frac{J(\bm{w})}{J^0} + \alpha g(\bm{w})^2 + \lambda g(\bm{w})
    \label{eq:lossFunction}
\end{equation}
where the parameters $\alpha$ and $\lambda$ are updated during each iteration, making the enforcement of the constraint stricter as the optimization progresses (see discussion below). Thus the overall framework is illustrated in \cref{fig:algoFlowChart}.

\subsection{Sensitivity Analysis}
\label{sec:sensAnalysis}
For minimizing the loss function in \cref{eq:lossFunction} we use  L-BFGS  \cite{liu1989limited}, a  well-known optimization technique. Since L-BFGS is a gradient-based optimizer, it requires sensitivities, i.e., derivative of the loss function in \cref{eq:lossFunction} with respect to the design variables  $\bm{w}$. Fortunately, one can exploit modern automatic differentiation frameworks \cite{chandrasekhar2021auto} to avoid manual sensitivity calculations. In particular, we use PyTorch \cite{NEURIPS2019_9015} as an implementation environment,  resulting in an end-to-end differentiable framework.

\subsection{Algorithm}
\label{sec:method_algo}
The  algorithm for the proposed framework is summarized below. First, we generate the dataset $\bm{D_m}$ consisting of permeability matrices (at various sizes) and contact area values for the pre-selected microstructures (line 2). Polynomials are  constructed to fit the data (line 3).

Next the domain is discretized for finite element analysis, and  the stiffness matrix components are computed (lines 4-5). The mesh is sampled at the center of each element (line 6); these serve as inputs to the NN. The augmented Lagrangian parameters $\alpha$, $\lambda$, the penalty parameter $p$ and NN weights $\bm{w}$ are  initialized (line 7).

In the main iteration, the design variables $\overline{\bm{\zeta}}$ are computed using the NN (line 9), followed by the computation of the permeability matrices for each element (line 10). These  are then used to construct the stiffness matrix and to solve for the velocity and pressure  (line 11 - line 12). Then the objective and contact area constraint are computed (lines 13 - 14), leading to the loss function (line 15). The  sensitivities are computed in an automated fashion (line 16). The weights $\bm{w}$ are then updated using  L-BFGS optimization scheme (line 17).  Finally the augmented Lagrangian multipliers and penalty parameters are updated (line 18 - 20). The process is repeated until termination, i.e., until the relative change in loss is below a certain threshold or the iterations exceed a maximum value. 

\begin{algorithm}[]
	\caption{GM-Flow}
	\label{alg:GM-Flow}
	\begin{algorithmic}[1]
		\Procedure{GMFlow}{$\Omega^0$ , BC,   $\Gamma^*$}
  
		\State $s \rightarrow \bm{D_m}  \quad s \in [0,1]$\Comment{calculate permeability and contact area at size instances} 
  
        \State $\bm{D_m} \rightarrow \bm{C}_m(s), {\Gamma}_m(s)$ \Comment{generate polynomial functions from dataset }\cref{fig:permVsSizeHalfFishscale}
        
		\State $\Omega^0 \rightarrow \Omega^0_h$ \Comment{discretize domain for FE} \label{alg:domainDiscretize}\cref{sec:Fluid_FEA}
		
		\State $\Omega^0_h \rightarrow \bm{A^\mu},\bm{A^\alpha},\bm{B}, \bm{a}$ \Comment{ compute stiffness matrices} \label{alg:stiffnessTemplates}\cref{eq:stiffness_term_defn}
		
		\State $\bm{x} = \{x_e,y_e\}_{e \in \Omega^0_h} \quad \bm{x} \in \mathbb{R}^{n_e \times 2}$ \Comment{elem centers; NN input} \label{alg:elemCenterComp}

		\State  epoch = 0; $\alpha = \alpha_0$; $p = p_0$; $\bm{w} = \bm{w}_0$ \Comment{initialization} \label{alg:initalizationParams}
		
		\Repeat \Comment{optimization (Training)}
		
		\State $NN(\bm{x} ; \bm{w}) \rightarrow \overline{\bm{\zeta}}(\bm{x})$ \Comment{fwd prop through NN} \label{alg:fwdPropNN}
		
		\State $\overline{\bm{\zeta}}(\bm{x}) \rightarrow \bm{C}_m(\bm{x}) $ \Comment{Permeability tensor  } \label{alg:PermeabilityTensorCompute}
  \cref{eq:effectiveCMatrix}
		
		\State$ \bm{C}_m(\bm{x}) \rightarrow \bm{K}, \bm{f} $ \Comment{ Stiffness matrix } \label{alg:stiffnessMatrix}
  \cref{eq:StiffnessMatrix}
				
		\State $\bm{K}, \bm{f}  \rightarrow \bm{S}$ \Comment{ solve \cref{eq:optimization_nn_govnEq}} \label{alg:feSolve}
		
		\State$\bm{K}, \bm{S} \rightarrow J$ \Comment{Objective, \cref{eq:optimization_nn_objective}} \label{alg:objectiveComp}

		\State $ \overline{\bm{\zeta}}, \Gamma^*  \rightarrow g_{\Gamma}$ \Comment{Contact area constraint }
  \cref{eq:effectivecontact area}\label{alg:perimCons}
		
		\State $J, g_{\Gamma} \rightarrow L$ \Comment{loss from \cref{eq:lossFunction}} \label{alg:lossCompute}
		
		\State $AD(L, \bm{w}) \rightarrow \nabla L $ \Comment{sensitivity analysis via Auto. Diff} \label{alg:autoDiff}
			 
		\State $\bm{w} , \nabla L \rightarrow \bm{w} $ \Comment{BFGS optimizer step}\label{alg:BFGSStep}
		
		\State $  \alpha + \Delta \alpha \rightarrow \alpha$ \Comment {increment penalty} \label{alg:OptPenaltyUpdate}
		
		\State $  \lambda + 2 \alpha g_\Gamma \rightarrow \lambda$ \Comment {increment Lagrange multiplier} 
  \cref{eq:lossFunction}\label{alg:OptLagrangiaUpdate}
		
		\State $   p + \Delta p \rightarrow p$ \Comment { continuation } \label{alg:materialPenaltyUpdate}
		
		\State $\text{epoch}++$
		
		\Until{ $|| \Delta L || < \Delta L_c^*$ or epoch < max\_epoch} \Comment{check for convergence}
		
		\EndProcedure
	\end{algorithmic}
\end{algorithm}


\section{Numerical Experiments}
\label{sec:experiments}

In this section, we conduct several experiments to illustrate the method and algorithm.  All experiments are conducted on a MacBook Air M2, using the PyTorch \cite{NEURIPS2019_9015} environment. The default settings are as follows:

\begin{itemize}
\setlength\itemsep{1em}
\item \textbf{Neural Network:} The NN comprises of 2 LeakyReLU-activated hidden layers with 20 neurons in each layer. This corresponds approximately to 4730 design variables. The initial values for $\bm{w}$ are determined via Xavier weight initialization \cite{glorot2010understanding}, with a seed value of 77.
\item \textbf{Candidate microstructures}: A set of 8 predefined microstructures (\cref{fig:microstructures}) is used in the experiments. Numerical homogenization is performed at six uniformly spaced sample points. A quintic polynomial is then used to interpolate the components of the sampled permeabiliity matrix.
\item \textbf{Material Penalization:} The penalization  $p$ is incremented every iteration by 0.02 using the continuation scheme \cite{sigmund1998numerical}, starting from a value of 1, with a maximum value of 8.
\item \textbf{Loss Function:} The initial constraint penalty is $\alpha_0$ = 0.05, and is increased by $\Delta \alpha = 0.15$ per epoch.
\item \textbf{Optimizer:} L-BFGS optimizer with a  strong Wolfe line search function is used \cite{nocedal1999numerical}. The maximum number of iterations (epochs) is set to 25. For convergence of optimization, we set change in loss $\Delta L_c^{*} = 10^{-5}$. 
 
\end{itemize}
For reference, the offline numerical homogenization of the 8 microstructures, at 6 different sizes, was performed in 16 seconds.

\subsection{Double-pipe Problem}
\label{sec:Double_pipe}
\subsubsection{Validation}
\label{sec:double_pipe_validation}
\begin{figure}[]
 	\centering
		\includegraphics[scale=0.45,trim={0 0 0 0}, angle= 270]{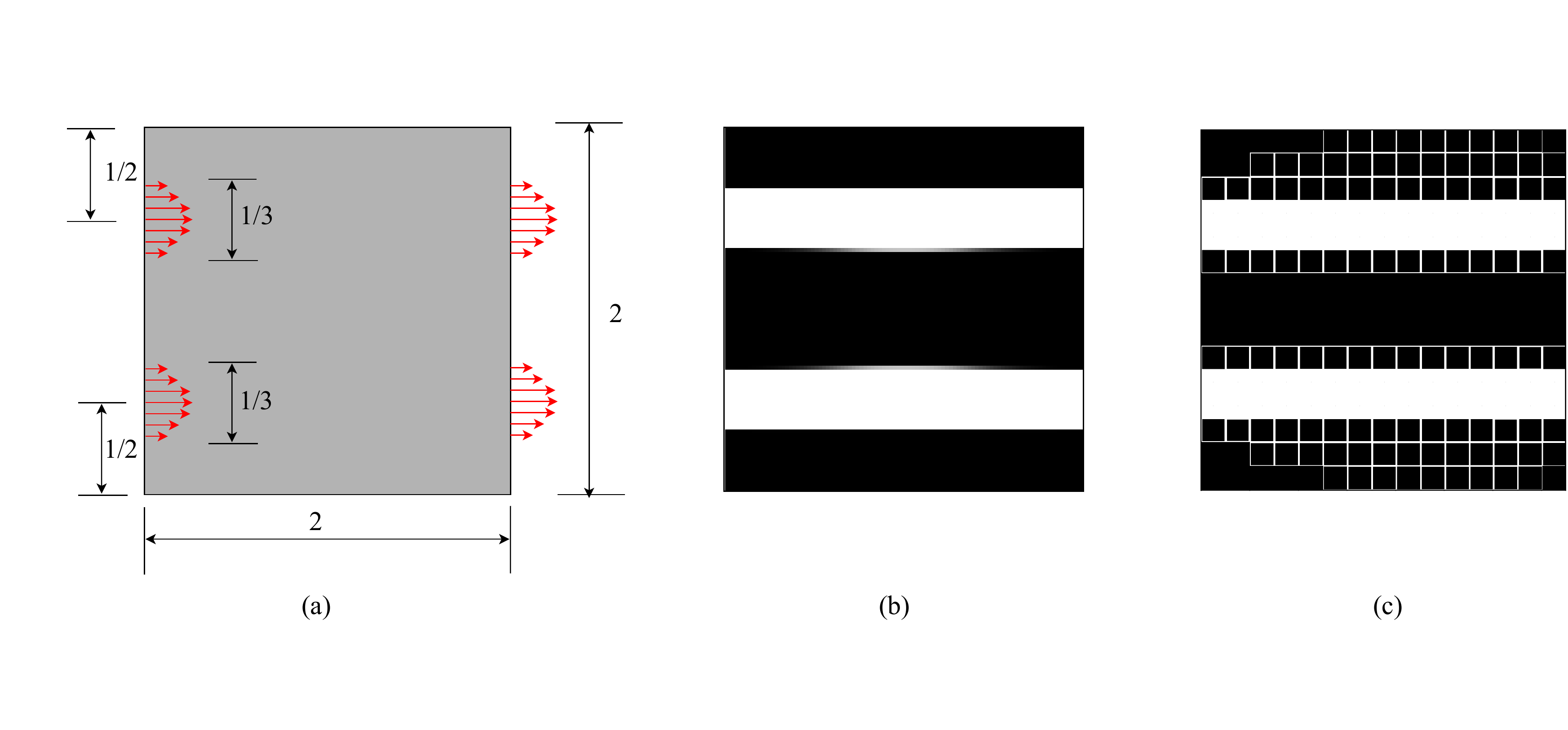}
 		\caption{Validation using the double pipe problem \cite{borrvall2003topology}.}
 		\label{fig:double_pipe_val}
 \end{figure}

In the first experiment, we validate the proposed method using the  double-pipe problem considered in \cite{borrvall2003topology}, and illustrated in \cref{fig:double_pipe_val}(a).  The objective is to find the optimal topology of $33\%$ fluid volume fraction that minimizes the dissipated power; contact area constraint is not imposed. The domain is discretized into 15x15 elements. The authors of \cite{borrvall2003topology} report the topology illustrated in \cref{fig:double_pipe_val}(b), with an objective value of $J=25.7$. In the proposed method, we use a single square microstructure (see \cref{fig:microstructures}), and arrive at the topology illustrated in \cref{fig:double_pipe_val}(c), with an objective of $J=27.4$. We observe partial infills in some of the cells since the size is not penalized towards $0,1$, i.e., we allow for intermediate sizes. 

\subsubsection{Impact of microstructures}
\label{sec:mstr_var_double_pipe}

Next we replace the square microstructure with three other microstructures, namely the circle, fish-body-2, and Mucosa-10 (see \cref{fig:microstructures}), one at a time, and study the impact of the microstructure on the dissipated power; the desired contact area is kept constant at $\Gamma^* = 60$. However, no constraint is imposed on the volume fraction. 

\begin{figure}[]
 	\centering
		\includegraphics[scale=0.45,trim={0 0 0 0}, angle= 270]{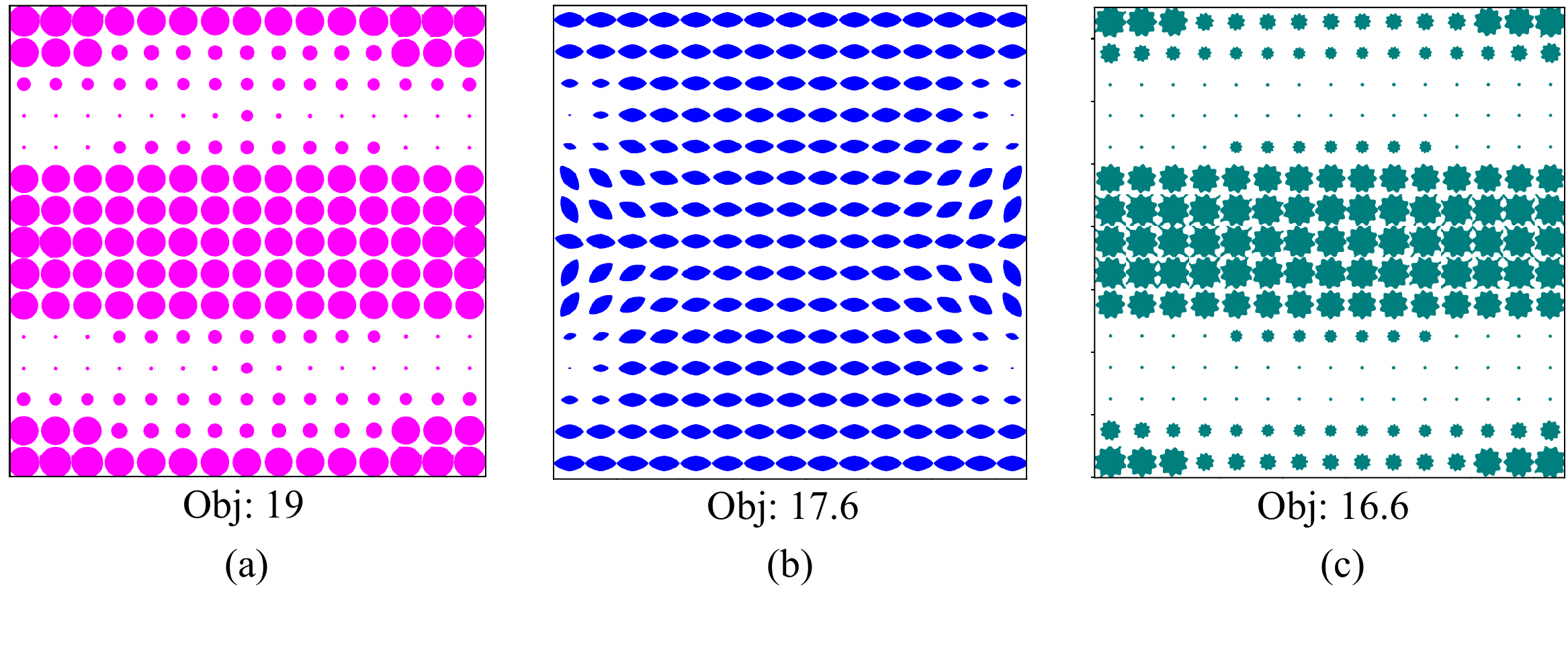}
 		\caption{Impact of microstructure on dissipated power with $\Gamma^* = 60$.}
 		\label{fig:double_pipe_single_mstr}
 \end{figure}

The resulting topologies are illustrated in \cref{fig:double_pipe_single_mstr}. Note that for the same contact area, the permeability of the fish-body-2 is higher than that of circle. This results in the fish-body-2 having a smaller dissipated power compared to the circle. On the other hand, despite a lower-permeability, Mucosa-10 performs better than the fish-body-2 since its contact area is significantly higher. 

\begin{figure}[]
 	\centering
		\includegraphics[scale=0.45,trim={0 0 0 0}, angle= 270]{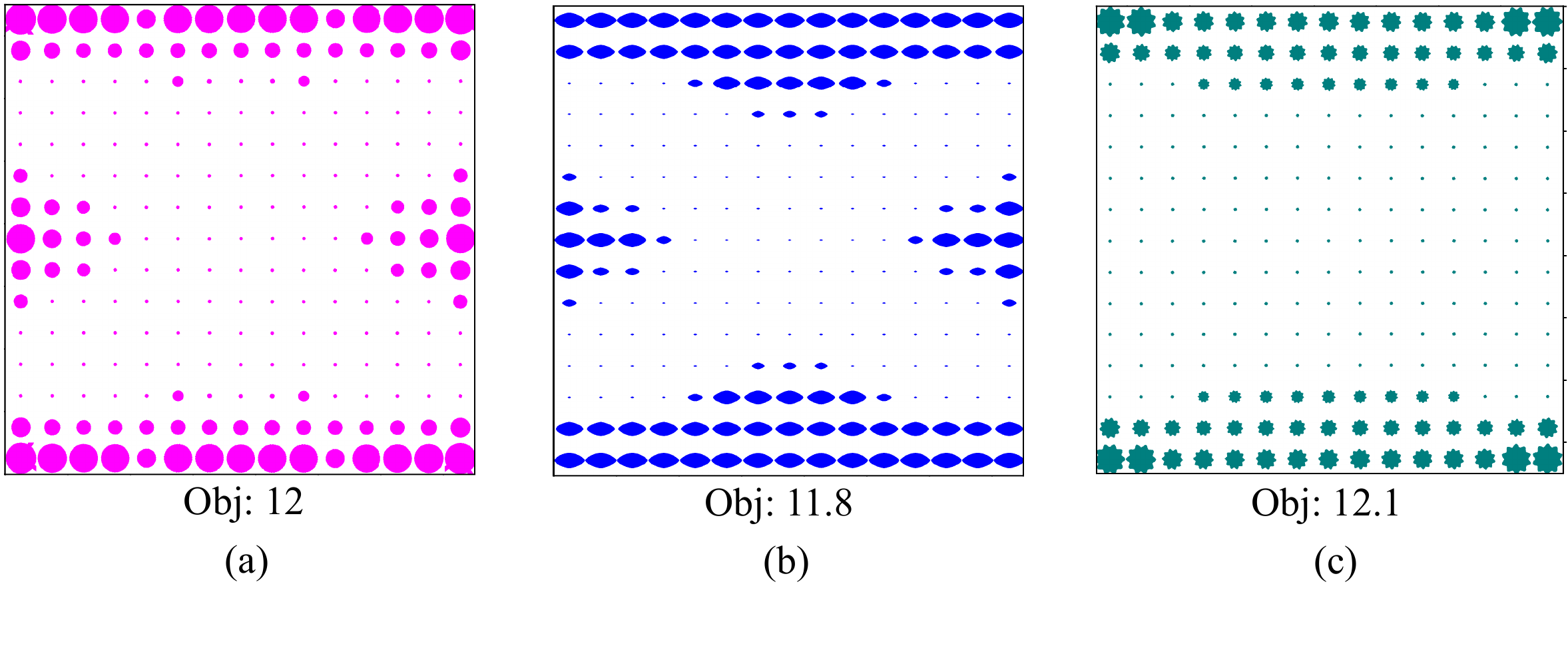}
 		\caption{Impact of microstructure on dissipated power with $\Gamma^* = 30$.}
 		\label{fig:double_pipe_single_mstr_30}
 \end{figure}

We then reduced the surface constraint to $\Gamma^* = 30$; the results are illustrated in \cref{fig:double_pipe_single_mstr_30}. Now, fish-body-2 and circle perform better than Mucosa-10. This experiment illustrates that the choice of a microstructure critically depends on the constraints imposed.

\subsection{ Diffuser Problem}
\label{sec:Validation}
\subsubsection{Validation}
\label{sec:diffuser_validation}
\begin{figure}[H]
 	\begin{center}
		\includegraphics[scale=0.45,trim={0 0 0 0}, angle= 270]{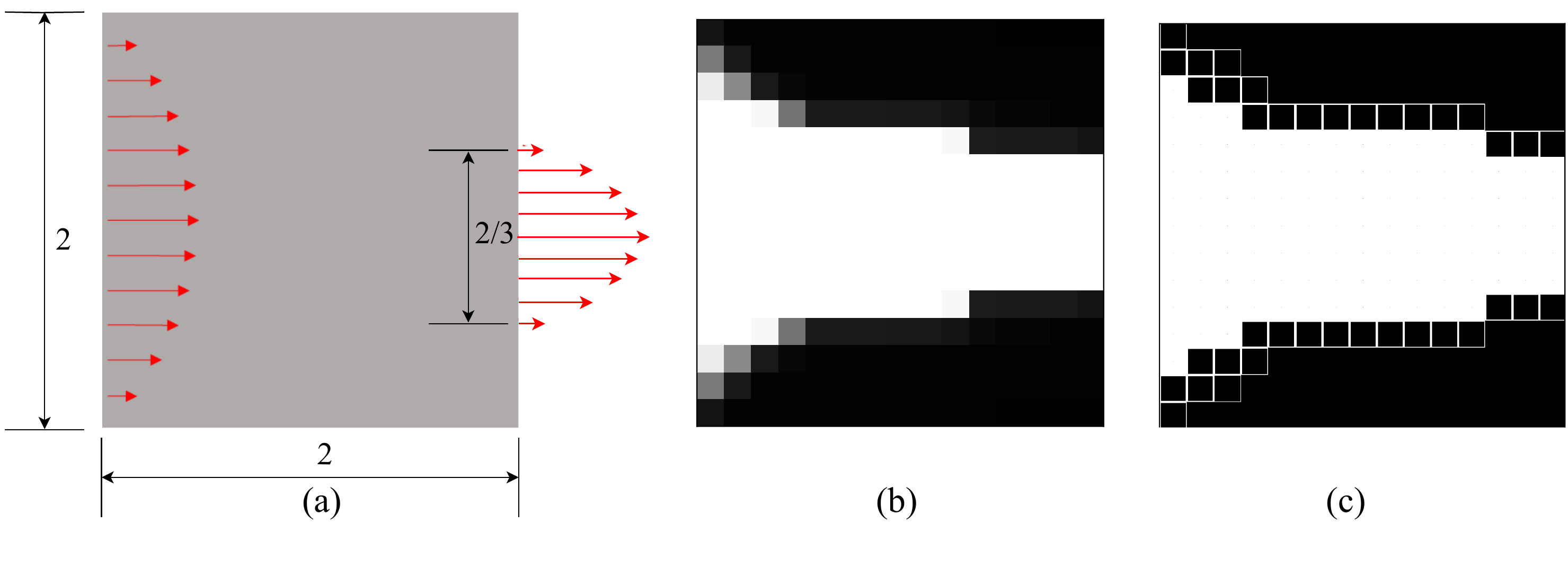}
 		\caption{(a) Diffuser problem. (b) Topology reported in \cite{borrvall2003topology}, (c) Topology generated via proposed method.}
 		\label{fig:0-1_val}
	\end{center}
 \end{figure}

For the second set of experiments, we consider the diffuser problem discussed in \cite{borrvall2003topology}, and illustrated in \cref{fig:0-1_val}(a).  In the first experiment, the objective is to find the optimal topology of $50\%$ fluid volume fraction that minimizes the dissipated power; contact area constraint is not imposed. The domain is discretized into 15x15 elements. The authors of \cite{borrvall2003topology} report the topology illustrated in \cref{fig:0-1_val}(b), with an objective value of $J=30.6$. Next, in the proposed method,  the square microstructure is once again used for optimization. The final topology is illustrated in \cref{fig:0-1_val}(c), consistent with \cref{fig:0-1_val}(b), with an objective of  $33.4$ \cite{borrvall2003topology}. The single scale optimization was performed in five seconds. 
 
 \subsubsection{Multiple microstructures}
\label{sec:No_microstructures}
\begin{figure}[]
 	\centering
		\includegraphics[scale=0.6,trim={0 0 0 0},clip]{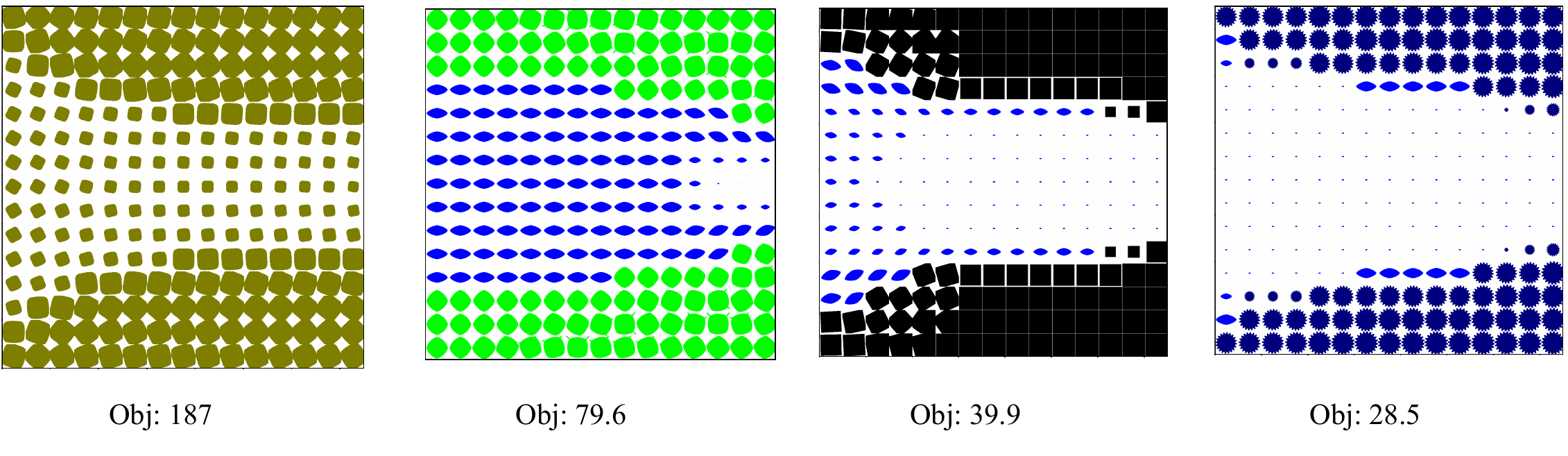}
 		\caption{Designs with 1, 3, 5 and 8 microstructures respectively.}
 		\label{fig:diffuser_mstrs_var}
 \end{figure}

Two central hypotheses of the current work is that one can achieve better designs with larger number of candidate microstructures, and that the framework is computationally insensitive to the number of candidates. To validate, we consider again the problem in \cref{fig:0-1_val}(a), but instead of the volume constraint, we impose a  contact area constraint of $\Gamma^* = 70$.  We consider the first $m$ microstructures, where $m = 1, 3, 5, 8$, and study the impact on the dissipated power and computational time. The resulting topologies are illustrated in \cref{fig:diffuser_mstrs_var}, as expected, the objective improves as we allow for larger number of microstructures. The computational time was approximately 62 seconds, independent of the number of microstructures.

 \subsection{ Bent Pipe Problem}
 \label{sec:bent_pipe}
 \subsubsection{Validation}
 \label{sec:val_bent_pipe}
\begin{figure}[]
 	\begin{center}
		\includegraphics[scale=0.65,trim={0 0 0 0},clip]{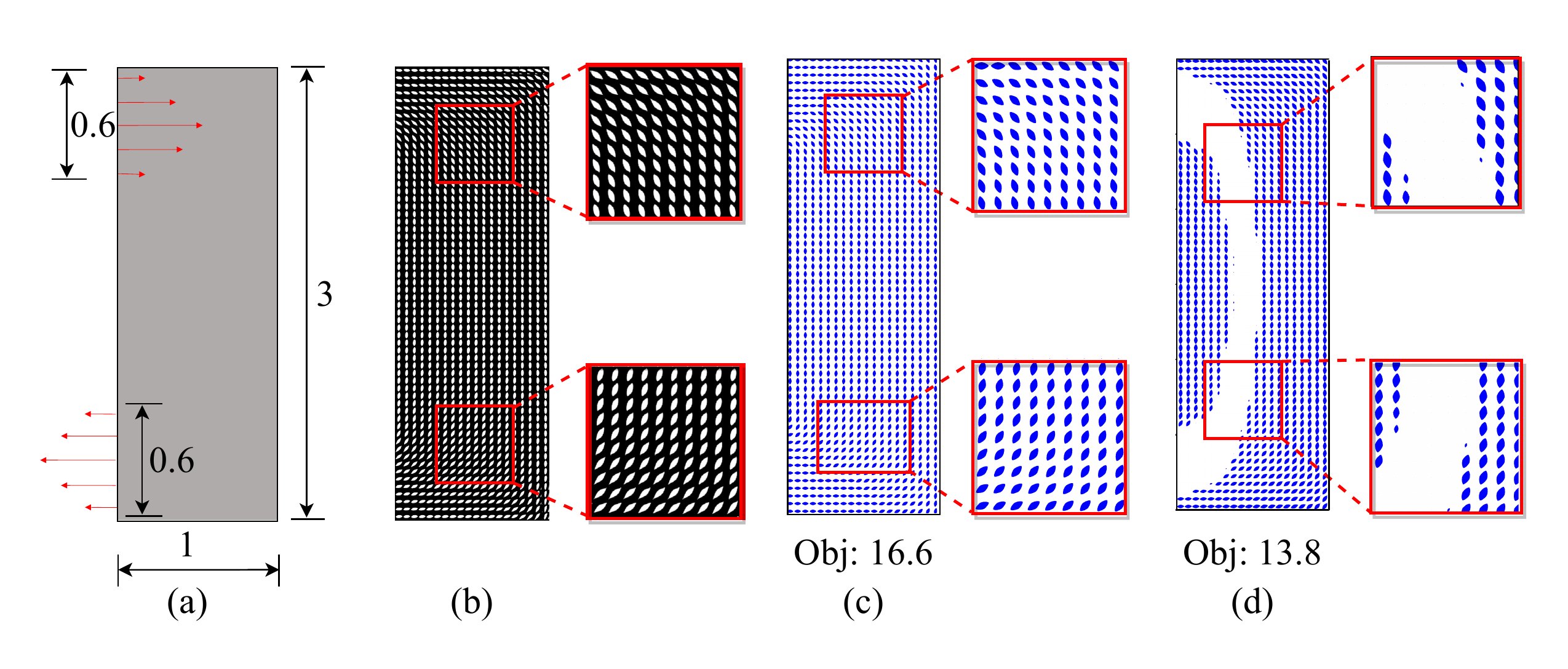}
 		\caption{Size fixed geometry validation (a) Domain and velocity boundary conditions (b) Solution reported in \cite{wu2019topology} (c) Topology generated via proposed method (d) solution via proposed method without rotation of microstructures.}
 		\label{fig:fix_size_val}
	\end{center}
 \end{figure}

Next, we consider the bent-pipe problem proposed in \cite{wu2019topology}, and illustrated in \cref{fig:fix_size_val}(a). The domain is discretized into 20x60 elements. In \cite{wu2019topology}, a two-scale topology optimization was carried out to minimize the dissipated power, with a constraint that the optimal microstructure must occupy exactly 25 percent of each unit cell. The reported topology is illustrated in \cref{fig:fix_size_val}(b); the final dissipated power was not reported. However, as noted in \cite{wu2019topology}, the computed microstructures resemble the fish-body. 

In the proposed method, the fish-body-2 microstructure was chosen  a priori, and its size was fixed  to occupy 25 percent of each unit cell, as in \cite{wu2019topology}. The orientation of each microstructure was optimized, resulting in the design  illustrated in  \cref{fig:fix_size_val}(c)  with the final dissipated power of  16.6. Not surprisingly, the final topology is similar to \cref{fig:fix_size_val}(b).

To improve on this design, we removed the constraint of the 25 percent unit-cell volume occupation, and instead imposed a total (global) volume constraint of 25 percent, i.e., we allowed the size of each microstructure to vary as well. The resulting topology is illustrated in \cref{fig:fix_size_val}(d). As one can observe, the dissipated power further reduces to 13.8.
\subsubsection{Fluid flow validation}
\label{sec:fluid_bent_pipe}
Next, for the above problem, we compute the pressure predicted using the GMTO framework and compare it against full-scale fluid flow simulation using  Ansys. Due to challenges in importing the geometry into ANSYS, the domain was discretized using a coarser mesh of 8x24 elements.

\begin{figure}[]
 	\begin{center}
		\includegraphics[scale=0.5,trim={0 0 0 0}]{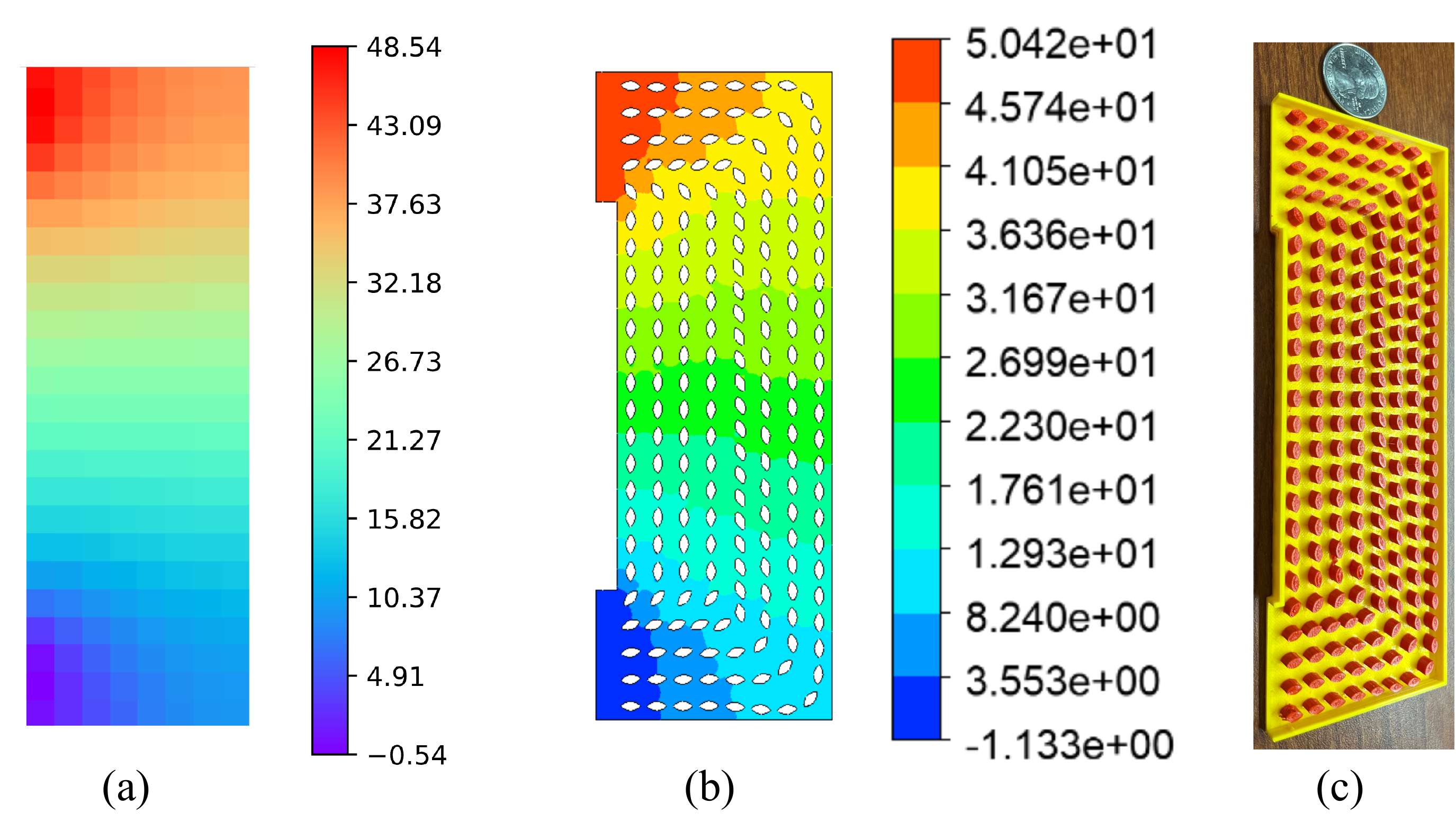}
 		\caption{Bent pipe pressure drop using: (a)  Proposed homogenization approach (b) Ansys (c) 3D printed design with microstructures.}
 		\label{fig:3d_print}
	\end{center}
 \end{figure}

The pressure prediction using the homogenization-based GMTO framework is illustrated in \cref{fig:3d_print}(a); the total pressure drop is approximately 49.0 Pascals. We then exported the optimized topology as an ".stl" file and imported it into ANSYS \cite{desalvo1979ansys} for full-scale fluid-flow simulation. The pressure drop predicted using ANSYS is illustrated in \cref{fig:3d_print}(b);  the total pressure drop is approximately 51.17 Pascals. The 3d printed part is illustrated in \cref{fig:3d_print}(c).
 
\subsubsection{Pareto Designs}
\label{sec:Pareto}

\begin{figure}[]
 	\begin{center}
		\includegraphics[scale=0.5,trim={0 0 0 0},clip]{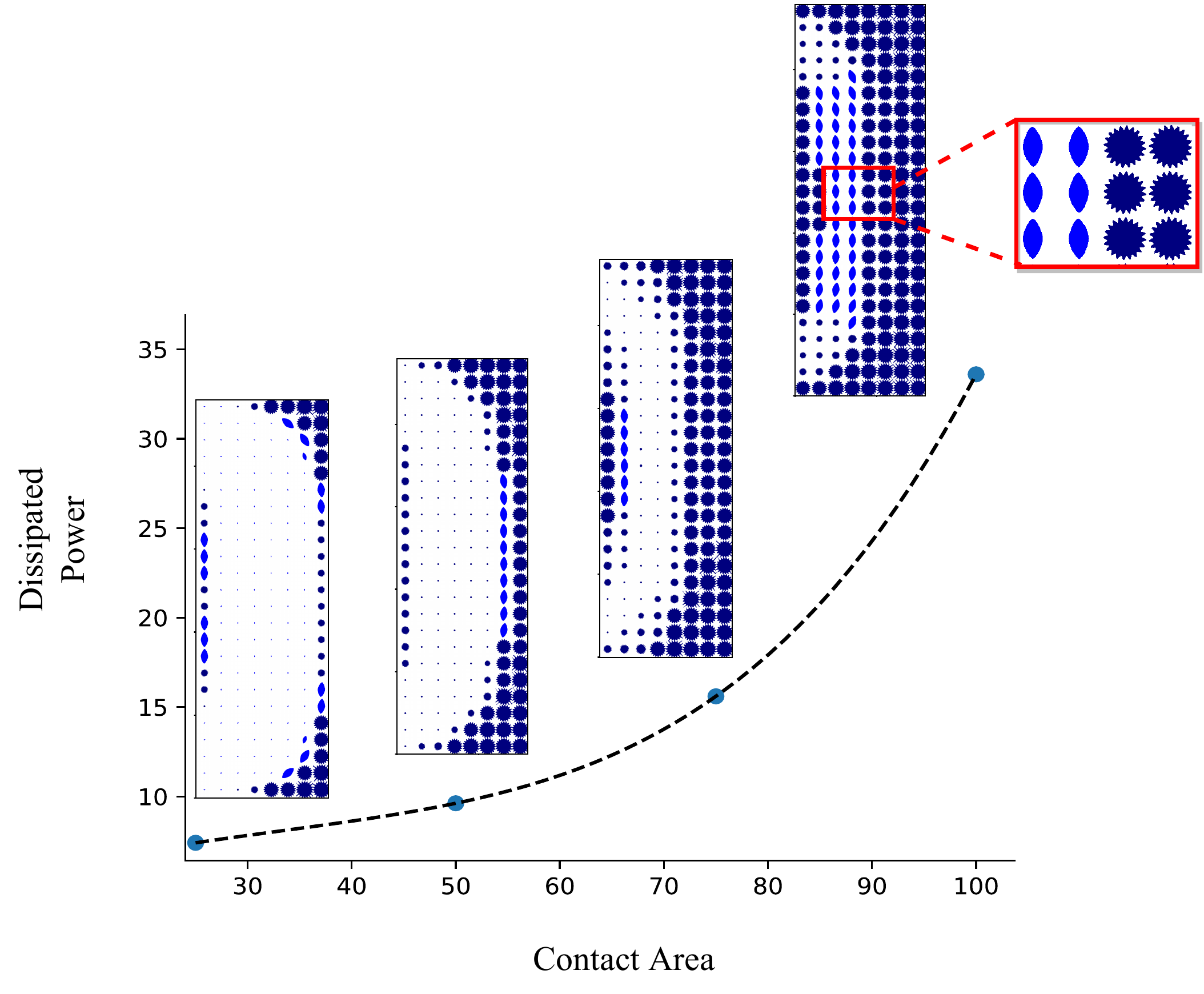}
 		\caption{Dissipated power versus contact area.}
 		\label{fig:ObjVsperim}
	\end{center}
 \end{figure}

Exploring the Pareto-front is critical in making design choices and understanding the trade-off between the objective (dissipated power) and constraint (contact area). To illustrate this, we computed the optimal topologies for various values of contact areas for the bent-pipe problem, by considering all 8 microstructures; the results are illustrated in \cref{fig:ObjVsperim}. Observe that as expected, the dissipated power increases with the contact area. Further, one can observe that regions with high fluid flow (see inset) are dominated by fish-body-2 (that exhibits high permeability), whereas regions with low fluid flow are dominated by  Mucosa-20 (that exhibits high contact area).

\subsubsection{Resampling }
\label{sec:Smooth}
\begin{figure}[H]
 	\begin{center}
		\includegraphics[scale=0.6,trim={0 0 0 0},clip]{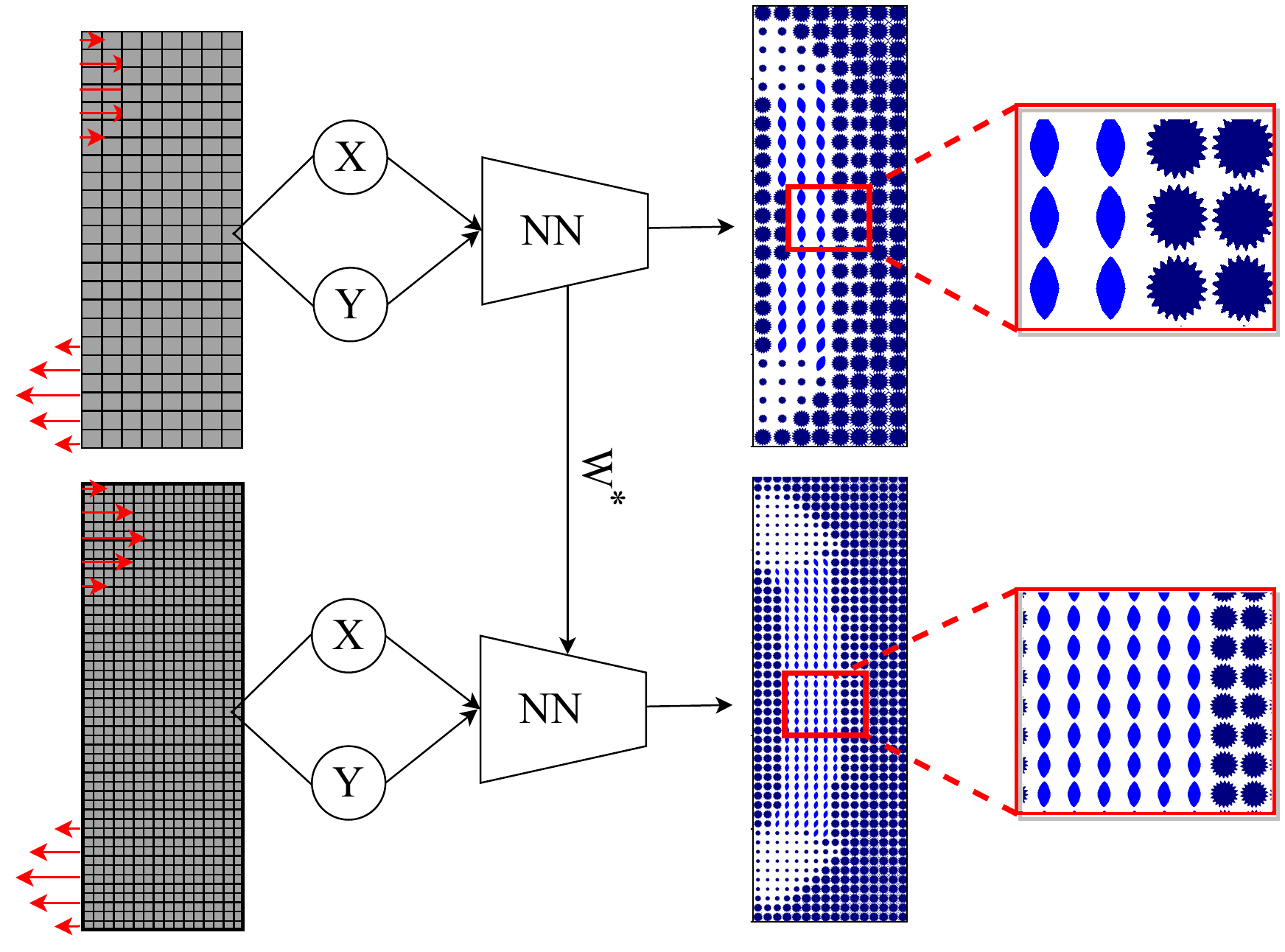}
 		\caption{Extraction of high resolution design through resampling. }
 		\label{fig:high res bent}
	\end{center}
 \end{figure}

A subtle but important aspect of the proposed  framework is that since the design fields (density, orientation, and size) are represented globally using the neural-network, one can obtain high-resolution topologies with no additional cost. To illustrate, suppose we have computed the optimal weights $\bm{w^*}$ and optimal topology for a given mesh discretization. We can then sample the domain at a higher resolution using the optimized weights $\bm{w^*}$, resulting in a more-detailed topology (Note that this is not a simple linear interpolation.) This is illustrated in \cref{fig:high res bent} where we optimize using a 8 × 24 mesh, and then re-sample using a 16 × 48 mesh.


\section{Conclusion}
\label{sec:conclusion}
A  graded multi-scale fluid flow topology optimization framework was proposed where homogenization  is  performed  off-line, followed by global optimization. Two-scale designs with high contact area and low dissipated power were generated through the framework. Furthermore, the  computational cost was found to be independent of the number of pre-selected microstructures. The neural-network configuration ensured that the partition of unity constraint was automatically satisfied. Finally, the PyTorch environment allowed for automated sensitivity computation. 

Currently, the framework is limited to microstructures with a single-size parameter. A contact area constraint, whose value was arbitrarily chosen, was imposed in the current work.  We plan to explore extensions to multi-physics problems such as convection-driven heat transfer problems \cite{dilgen2018density} where the contact area constraint is determined through the underlying physics. Furthermore, it will be interesting to combine the proposed framework with data-driven methods \cite{wang2022data}.

\section*{Compliance with ethical standards}
\label{sec:ethics}
The authors declare that they have no conflict of interest.
\section*{Replication of Results}
\label{sec:replic}
The Python code pertinent to this paper is available at \href{https://github.com/UW-ERSL/FluTO}{github.com/UW-ERSL/FluTO}.
\section*{Acknowledgement}
The authors would like to thank the support of National Science Foundation through grant CMMI 1561899. The authors acknowledge Subodh Subedi for helping with the 3D printing.

\bibliographystyle{plain}
\bibliography{references}  

\end{document}